\pdfoutput=1

\documentclass[draft]{article}
\usepackage[T1]{fontenc}
\usepackage[UKenglish]{babel}

\usepackage{Factory, Math, Classic, Theorema, Styl}

\newcommand{\ieme}{$^{\text{e}}$}
\newcommand{\no}{n$^{\text{o}}$}

\usepackage[all, 2cell]{xy} \UseAllTwocells \newcommand{\dotar}{\ar@{.>}}

\newcommand{\Sets}{\mathfrak{Set}}		
\newcommand{\Mod}{\mathfrak{Mod}}			
\newcommand{\SMod}{{}_S\Mod} 			
\newcommand{\CAlg}{\mathfrak{CAlg}}

\newcommand{\NAlg}{\mathfrak{BAlg}}
\newcommand{\freenum}[2][\Z]{#1\binom{#2}{-}}

\newcommand{\B}{\mathbf{B}}

\newcommand{\de}{\diamond} 			
\DeclareMathOperator*{\De}{ \lozenge } 
\newcommand{\dev}[3]{#1(#2\de\cdots\de #3)}   	
\newcommand{\Dev}[2]{#1\left(\De #2\right)}

\usepackage{fancyhdr}
\begin{document}

\titul{POLYNOMIAL MAPS \\ OF MODULES} 
\auctor{Qimh Richey Xantcha%
\thanks{\textsc{Qimh Richey Xantcha}, Uppsala University: \texttt{qimh@math.uu.se}}}
\datum{\today}
\maketitle

\bigskip 

\epigraph{\begin{vverse}[And God said unto the animals: ``Be fruitful and multiply.''] 
And God said unto the animals: ``Be fruitful and multiply.'' \\[0.5ex]

But the snake answered: ``How could I? I am an adder!''\footnote{In some retellings of this myth, 
it is further reported that God constructed a wooden table for the snakes to crawl upon, since even 
adders can multiply on a log table. 
God is not assumed to be familiar with tensor products.}
\end{vverse}
}

\bigskip 

\begin{argument} \noindent
The article focuses on three different notions of polynomiality for maps of modules. 
In addition to the polynomial maps studied by Eilenberg and Mac Lane and the strict polynomial maps 
(``lois polynomes'') considered by Roby, 
we introduce \emph{numerical maps} and investigate their properties. 

Even though our notion require the existence of binomial co-efficients in the base ring, we argue that 
it constitutes the correct way to extend
Eilenberg and Mac Lane's polynomial maps of abelian groups to incorporate modules 
over more general rings. 
The main theorem propounds that our maps admit a description corresponding, word by word,  
to Roby's definition of strict polynomial maps. 

\end{argument}

\bigskip 

\noindent 
Linear maps have been extensively studied in Modern Algebra; polynomial maps less so. 
This was the rationale behind Eilenberg and Mac Lane's introduction in \cite{EM}  of 
\emph{polynomial maps} of abelian groups. Passi (\cite{PassiI}, \cite{PassiII}, \cite{PassiBook}) 
and Buckley (\cite{Buckley}) have continued this line of enquiry, 
expanding the concept to accommodate maps from possibly non-abelian groups (or even monoids) 
to abelian groups. Passi's excellent survey \cite{PassiSurvey} may serve as an introduction to the subject. 

While the approach taken, regulating the ``non-additivity'' 
of the map under consideration, will remain valid 	
for modules over any ring, it will clearly be deficient, for the simple and obvious reason that it does not take scalar multiplication into account. 

Roby was then led to consider \emph{strict polynomial maps} of modules 
over an arbitrary commutative base ring $\B$.
According to \cite{Roby}, a strict polynomial map of $\B$-modules is a natural transformation 
$$
\phi\colon M\otimes_\B - \to N\otimes_\B -
$$ between functors ${}_\B\CAlg\to \Sets$, 
where ${}_\B\CAlg$ denotes the category of 
commutative, unital algebras over $\B$. 
A broad and all-encompassing notion of polynomial maps as types of natural 
transformations (Definition \ref{Def: Polynomiality}) can be extracted 
as the essence of this.

One observes that strict polynomial maps are not determined by the underlying set maps. 
The simplest example is probably the following. Let $\B=\Z$, and define 
$$ 
\phi_A\colon \Z/2\otimes A \to \Z/2\otimes A, \qquad 1\otimes x\mapsto 1\otimes x(x-1). 
$$
This is a non-trivial strict polynomial map of degree $2$, and its underlying map is zero. 
We see here at play the well-known distinction between polynomials and polynomial maps, 
the former class being richer than the latter. 
The point is that the strict polynomial structure provides extra data, which make the zero map strict polynomial of degree 2 in a non-trivial way.

In this paper, we introduce  
\emph{numerical maps}, which we believe furnish the proper 
way to extend Eilenberg and Mac Lane's weak notion of polynomiality to more general base rings, handling, as it 
does, scalar multiplication in a natural way. The key point is   
that the base ring $\B$ is required to possess binomial co-efficients. That is, it should be a 
\emph{binomial ring} in the sense of Hall~\cite{Hall}. 

Let $\phi\colon M\to N$ be a module map, and let $\phi(x_1\de \cdots \de x_{n+1})$ denote the deviation 
in $n+1$ variables (see Definition \ref{D: Deviations} below). 
The vanishing of the latter is polynomiality in the sense of Eilenberg \& Mac Lane, to which we 
adjoin an auxiliary condition:

\begin{intdefinition}[\ref{Def: Numerical Maps}] 
The map $\phi\colon M\to N$ is \textbf{numerical} of degree (at most) $n$ if it satisfies the following two equations: 
\begin{gather*}
\dev{\phi}{x_1}{x_{n+1}} =0, \qquad x_1,\dots,x_{n+1}\in M \\
\phi(rx) = \sum_{k=0}^n \binom{r}{k}\Dev{\phi}{_k x}, \qquad r\in\B,\ x\in M.
\end{gather*} 
\end{intdefinition}

It is easy to prove that, over the integers, the second condition is implied by the first, 
so that the concepts of polynomial and numerical map co-incide. 

Likewise, it follows readily that a strict polynomial map is also numerical of the same degree 
(provided of course the base ring be binomial). 
If the base ring is a $\Q$-algebra, then, comparing Definition \ref{D: SPol} with Theorem \ref{Th: Num = NAlg} below, 
the two concepts will be seen to co-incide, every algebra being binomial in this case. 

We exhibit, in Theorem \ref{Th: Universal}, the universal numerical map of each degree, 
and further prove the following pleasant formula: 

\begin{inttheorem}[\ref{Th: Poly Map Eq 2}] 
The map $\phi\colon M\to N$ is numerical of degree $n$ if and only if  
$$ 
\phi(a_1x_1\de\cdots\de a_kx_k) = 
\sum_{\substack{ m_1+\cdots+m_k\leq n \\ m_1,\dots,m_k\geq 1}} \binom{a_1}{m_1}\cdots\binom{a_k}{m_k} 
\phi\left(\De_{m_1} x_1 \de\cdots\de \De_{m_k} x_k\right) 
$$ 
for any $a_i\in\B$ and $x_i\in M$.
\end{inttheorem}

The next theorem should be compared with Theorem \ref{S: Roby}, 
which gives the corresponding property for strict polynomial maps. Buckley \cite{Buckley} provided the first such 
characterisation of polynomial maps.

\begin{inttheorem}[\ref{Th: Char of Poly Maps}] 
Let $\phi\colon M\to N$ be numerical of degree $n$ and  
let $u_1, \ldots, u_k \in M$. For any natural numbers $m_1,\dots,m_k$ there is 
an element $v_{(m_1,\dots,m_k)}\in N$, non-zero only if $m_1+\cdots+m_k\leq n$, such that  
\beq 		\label{E: Char of Poly Maps}
\phi(a_1u_1+\cdots+ a_ku_k)= \sum \binom{a_1}{m_1}\cdots \binom{a_k}{m_k} v_{(m_1,\dots,m_k)}
\eeq
for any $a_1,\dots,a_k\in \B$. The elements $v_{(m_1,\dots,m_k)}$ are uniquely defined and, explicitly, 
$$
v_{(m_1,\dots,m_k)} = \phi\left( \De_{m_1} u_1 \de \cdots \de \De_{m_k} u_k\right).
$$

Conversely, the existence of elements $v_{(m_1,\dots,m_k)}$ with the above property 
ensures that $\phi$ is numerical of degree $n$.
\end{inttheorem}


Remarkably, these numerical maps will be seen to fit perfectly
into the frame-work erected by Roby for strict polynomial maps.  
Our final theorem will provide an instant unification of the two notions: 

\begin{inttheorem}[\ref{Th: Num = NAlg}] 
The map $\phi\colon M\to N$ is numerical of degree $n$ if and only if it can be extended to a (unique)
natural transformation
$$
\phi\colon M\otimes_\B - \to N\otimes_\B -,
$$
of degree $n$, between functors ${}_\B\NAlg\to \Sets$, where ${}_\B\NAlg$ denotes the 
category of binomial algebras over the ring $\B$.
\end{inttheorem}

One point may deserve further elaboration, particularly in light of the recent attempt \cite{Hartl} 
by Gaudier \& Hartl to 
define quadratic maps for arbitrary modules. 
In the original case considered by Eilenberg \& Mac Lane, viz.~$\B=\Z$, 
a polynomial map $\phi\colon M\to N$ of abelian groups
is equivalent to a natural transformation 
$$
\phi\colon M\otimes_\Z - \to N\otimes_\Z -
$$
of functors with source ${}_\Z\NAlg$, the category of binomial rings. 
Thus, even in the case of abelian groups, where no binomial co-efficients \emph{a priori} appear, 
binomial rings will none the less enter, in a canonical fashion, to render the definition 
akin to Roby's. 

This research was carried out at Stockholm University under the eminent supervision of
Prof.~Torsten Ekedahl. 
We would also like to thank Dr Christine Vespa for innumerous and invaluable comments 
on the manuscript.

\setcounter{section}{-1} 
\section{Preliminaries}  			\label{Sec: Preliminaries}

For the entirety of this article, $\B$ shall denote a fixed base ring of scalars, assumed commutative and unital. 
All modules, homomorphisms, and tensor products shall be taken over this $\B$, unless otherwise stated. 
We let $\Mod={}_\B \Mod$ denote the category of (unital) modules over this ring.

A \emph{map} of modules shall always denote an \emph{arbitrary} map --- in general non-linear. 
By the term ``homomorphism'', we shall, of course, mean a linear map. 

The deviations of a map were first defined by Eilenberg \& Mac Lane \cite{EM}. 
Let $[n]$ denote the set $\{1,\dots,n\}$. 

\bdf \label{D: Deviations}
Let $\phi\colon M\to N$ be a map of modules. The $n$'th \textbf{deviation} of $\phi$ is the map 
$$ \dev{\phi}{x_1}{x_{n+1}} = \sum_{I\subseteq[n+1]} (-1)^{n+1-\abs{I}} \phi\left(\sum_{i\in I} x_i\right) $$
of $n+1$ variables. 
\edf

As a special case, we have $\phi(\de )=\phi(0)$. We agree to abbreviate
$$
\Dev{\phi}{_n x} = \phi(\underbrace{x\de\cdots\de x}_{n}).
$$

The following formula is an immediate consequence of the definition, and will be found very useful in what follows.  

\begin{lemma*}			\label{L: Dev Lemma}
$$ 
\phi(x_1+\dots+ x_{n+1}) = \sum_{I\subseteq[n+1]} \Dev{\phi}{_{i\in I} x_i}.  
$$
\end{lemma*}


\bdf
The map $\phi\colon M\to N$ is \textbf{polynomial} of degree (at most) $n$ if its $n$'th deviation vanishes: 
$$
\dev{\phi}{x_1}{x_{n+1}} =0 
$$
for any $x_1,\dots,x_{n+1}\in M$.
\edf


We now present Roby's notion of strict polynomial map, as given in section 1.2 of \cite{Roby}. 
Let $\CAlg={}_\B\CAlg$ denote the category of commutative, unital algebras over the base ring $\B$, 
and let $\Sets$ denote the category of sets.  

\bdf 		\label{D: SPol}
A natural transformation 
$$\phi\colon M\otimes - \to N\otimes -$$ between functors $\CAlg\to \Sets$, 
will be called a \textbf{strict polynomial map}. 
\edf


\begin{reftheorem}[\cite{Roby}, Théorème \textsc{i}.1]		\label{S: Roby}
Let $\phi\colon M\to N$ be a strict polynomial map and let $u_1,\ldots,u_k\in M$. 
For any natural numbers $m_1,\dots,m_k$ there is an element 
$v_{(m_1,\dots,m_k)}\in N$ (only finitely many of these being non-zero)  
such that 
$$ 
\phi(u_1\otimes a_1 +\cdots +u_k\otimes a_k) 
= \sum v_{(m_1,\dots,m_k)}\otimes a_1^{m_1}\dots a_k^{m_k}, 
$$
for all $a_j$ in all (commutative, unital) algebras. Moreover, the elements $v_{(m_1,\dots,m_k)}$ are unique.
\end{reftheorem}



\section{Numerical Maps} 

The base ring $\B$ of scalars will now be assumed \emph{binomial} in the sense of \cite{Hall}, 
by which is 
meant a commutative ring with unity, equipped with binomial co-efficients. 
One usually 
postulates $\B$ to be torsion-free and closed in $\Q\otimes_\Z \B$ under the operations 
$$ 
r\mapsto \frac{r(r-1)\cdots(r-n+1)}{n!}. 
$$
The theory of these rings has been expounded in \cite{Binomial} and \cite{BR}.

\bdf \label{Def: Numerical Maps} 
Let $M$ and $N$ be modules over $\B$.
The map $\phi\colon M\to N$ is \textbf{numerical} of degree (at most) $n$ if it satisfies the following two equations: 
\begin{gather*}
\dev{\phi}{x_1}{x_{n+1}} =0, \qquad x_1,\dots,x_{n+1}\in M \\
\phi(rx) = \sum_{k=0}^n \binom{r}{k}\Dev{\phi}{_k x}, \qquad r\in\B,\ x\in M.
\end{gather*} 
\edf

\bex
When $n=0$ the above equations read: 
\begin{gather*}
\phi(x_1)-\phi(0)=\phi(\de x_1) =0, \qquad\qquad
\phi(rx) = \binom r0 \phi(\de) = \phi(0).
\end{gather*}
Consequently, the map $\phi$ is of degree $0$ if and only if it is constant.
\eex

\bex
When $n=1$, the equations read: 
\begin{gather*}
\phi(x_1+x_2)-\phi(x_1)-\phi(x_2)+\phi(0)=\phi(x_1 \de x_2) =0, \\
\phi(rx) = \binom r0 \phi(\de) + \binom r1 \phi(\de x) = \phi(0) + r(\phi(x)-\phi(0)).
\end{gather*}
The map $ \psi(x) = \phi(x) - \phi(0) $
is then a module homomorphism. Conversely, any translate of a module 
homomorphism is numerical of degree $1$. 
\eex

\bex
Let $\B=\Z$. It is a well-known fact, and not difficult to prove, that the numerical (polynomial) maps 
$\phi\colon\Z\to\Z$ of degree $n$ are of the form 
$\phi(x) = \sum_{k=0}^n c_k\binom{x}{k}$
(binomial polynomials of degree $n$). 
\eex

\bex 
Let $\B=\Z$. The map 
$$ \xi\colon \Z \otimes - \to \Z \otimes - , \qquad 1\otimes x\mapsto 1\otimes \binom x2, $$
is numerical of degree $2$, but not strict polynomial of any degree. The following tentative diagram, where $\beta\colon t\mapsto a$, indicates 
the impossibility of defining $\xi_{\Z[t]}$. 
$$ 
\xymatrix{
\Z\otimes \Z[t] \ar[r]^{\xi_{\Z[t]}} \ar[d]_{1\otimes\beta}  & \Z\otimes \Z[t] \ar[d]^{1\otimes\beta} & 1\otimes t \ar[d] \ar[r] & \text{?} \ar[d] \\
\Z\otimes \Z \ar[r]_{\xi_\Z} & \Z\otimes \Z & 1\otimes a \ar[r] & 1\otimes \binom a2  \\
} 
$$
Note that $\Z[t]$ is not a binomial ring, there being no such thing as $\binom t2$. 
\eex

%

\section{Polynomiality} 

Let us now conduct an investigation of \emph{polynomiality} at its most general, 
and indicate how this perspective provides a unifying view. 

Let $D\subseteq\Mod$ be a \emph{finitary algebraic category}, by which is simply meant an equational 
class in the sense of universal algebra 
(and hence a variety of algebras by the famous \textsc{hsp} Theorem). 
Since $D$ is a subcategory of $\Mod$, the objects of $D$ are first of all $\B$-modules, possibly equipped 
with some extra structure. We do not assume $\B$ is binomial in this section.

For a set of variables $V$, let $\gen{V}_D$ denote the free algebra on $V$ in $D$. 
That the free algebra exists is a basic fact of universal 
algebra; confer \cite{Universal}.

\bdf Let $M$ be a module, not necessarily in $D$. An element of  
$$
M\otimes \gen{  x_1,\ldots,x_k}_D
$$
is called a \textbf{$D$-polynomial} over $M$ in the variables $x_1,\ldots,x_k$.

A \textbf{linear form} over $M$ in these variables is a polynomial of the form $\sum u_j\otimes x_j$, 
for some $u_j\in M$. 
\edf

\bth[The Polynomiality Principle]
Let two modules $M$ and $N$ be given, and a family of maps 
$$\phi_A\colon M\otimes A \to N\otimes A, \quad A\in D.$$ 
The following statements are equivalent: \label{Th: Polynomiality}
\balph
\item For every $D$-polynomial $p(x)=p(x_1,\ldots,x_k)$ over $M$, there is a unique $D$\hyp polynomial 
$q(x)=q(x_1,\ldots,x_k)$ over $N$, such that for all $A\in D$ and all $a_j\in A$, 
$$
\phi_A(p(a)) = q(a) .
$$
\item For every linear form $l(x)$ over $M$, there is a unique $D$-polynomial 
$q(x)$ over $N$, such that for all $A\in D$ and all $a_j\in A$, 
$$
\phi_A\left(l(a)\right) = q(a) .
$$
\item The map $$\phi\colon M\otimes- \to N\otimes -$$ is a natural transformation of functors $D\to \Sets$. 
\ealph \hfill
\eth

\bpr
A trivially implies B. Suppose statement B holds, and consider an homomorphism $\chi\colon A\to B$, 
along with finitely many elements $u_j\in M$. Define 
$$l(x)=\sum u_j\otimes x_j,$$
and find the unique $D$-polynomial $q$ satisfying B. Then, for any $a_j\in A$, there is a commutative diagram of the following form, proving that $\phi$ is natural: 
$$ \xymatrix{
M\otimes A \ar[r]^{\phi_A} \ar[d]_{1\otimes\chi}  & N\otimes A \ar[d]^{1\otimes\chi} & \sum u_j\otimes a_j \ar[d] \ar[r] & q(a) \ar[d] \\
M\otimes B \ar[r]_{\phi_B} & N\otimes B & \sum u_j\otimes \chi(a_j) \ar[r] & q(\chi(a)) \\
} $$
Thus, condition C holds. 

Finally, suppose $\phi$ natural. We shall prove condition A. Given a $D$\hyp polynomial
$$p(x)\in M\otimes \gen{ x_1,\ldots,x_k}_D,$$ 
define 
$$
q(x)=\phi_{\gen{ x_1,\ldots,x_k}_D}\left(p(x)\right).
$$
For any $A\in D$ and $a_j\in A$, define the homomorphism 
$$\chi\colon \gen{ x_1,\ldots,x_k}_D\to A, \quad x_j\mapsto a_j.$$ 
Then since $\phi$ is natural, the following diagram commutes: 
$$ \xymatrix{
M\otimes \gen{  x_1,\ldots,x_k} \ar[rr]^{\phi_{\gen{ x_1,\ldots,x_k} }} \ar[d]_{1\otimes \chi}  & 
& N\otimes \gen{ x_1,\ldots,x_k}\ar[d]^{1\otimes\chi} & \ar[d] p(x) \ar[r] & q(x) \ar[d] \\
M\otimes A \ar[rr]_{\phi_A} && N\otimes A & p(a) \dotar[r] & q(a) 
} $$
The uniqueness of $q$ is evident, proving A. 
\epr

\bdf \label{Def: Polynomiality}
When the conditions of the theorem are fulfilled, we call $\phi$ a \textbf{$D$-polynomial map} from $M$ to $N$. 
\edf

According to part B of the theorem, $\phi_A$ maps
$$ \sum u_j\otimes a_j \mapsto q(a)  $$
for some (unique) $D$-polynomial $q$. In naïve language, the Polynomiality Principle  amounts to the following. 
\emph{If we want the co-efficients $a_j$ (in some algebra) of the module elements $u_j$ to transform according to certain operations, 
the correct setting is the category of algebras using these same operations.}

\bex
A \emph{$\Mod$-polynomial map} $\phi\colon M\to N$ is just a linear transformation $M\to N$. This is because, by B above, $\phi_\B$ will map 
$\sum u_j\otimes r_j$ to $\sum v_j\otimes r_j$ for all $r_j\in \B$, and such a map is easily seen to be linear. Conversely, any module homomorphism
induces a natural transformation $M\otimes -\to N\otimes -$. 
\eex

\bex
Let $S$ be a $\B$-algebra; then $\SMod\subseteq \Mod$. An \emph{$\SMod$\hyp polynomial map} 
$M\to N$ is a transformation 
$$ 
M\otimes A\to N\otimes A 
$$ 
which is natural in the $S$-module $A$. This is the same as a natural transformation 
$$
(M\otimes S) \otimes_S - \to (N\otimes S) \otimes_S -,
$$
which is an $\SMod$-polynomial map $M\otimes S\to N\otimes S$; or, as noted in the previous example, an $S$-linear map from $M\otimes S$ to $N\otimes S$.
\eex

\bex
The \emph{$\CAlg$-polynomial maps} are precisely the strict polynomial ones.  
Referring to Theorem \ref{S: Roby}, the equation  
$$
\phi(u_1\otimes a_1 +\cdots +u_k\otimes a_k) 
= \sum v_{(m_1,\dots,m_k)}\otimes a_1^{m_1}\dots a_k^{m_k}
$$
shows that, intuitively, the co-efficients of the elements $u_j$ ``transform as ordinary polynomials''.
\eex

Let now $\B$ be binomial, and consider the category $\NAlg={}_\B\NAlg$ of binomial 
algebras over $\B$.
According to the definition, a \textbf{$\NAlg$-polynomial map} is a natural transformation 
$$
\phi_A\colon M\otimes A \to N\otimes A, \quad A\in \NAlg.
$$ 
The Polynomiality Principle provides, for every linear form $\sum_{j=1}^k u_j\otimes x_j$ over $M$, 
a unique numerical polynomial 
$$
\sum v_{(m_1,\dots,m_k)}\otimes \binom{x_1}{m_1}\cdots\binom{x_k}{m_k}
$$ 
over $N$, with the property that, for all binomial algebras $A$ and all $a_j\in A$,
\beq 			\label{Eq: Num} 
\phi_A\left(u_1\otimes a_1 + \dots + u_k\otimes a_k\right) = 
\sum v_{(m_1,\dots,m_k)}\otimes \binom{a_1}{m_1}\cdots\binom{a_k}{m_k}.
\eeq
Intuitively, the co-efficients of the elements $u_j$ ``transform as numerical (binomial) polynomials''.

While the right-hand side of equation \eref{Eq: Num} is finite, there is no reason to expect 
a uniform upper bound for its degree. Roby refers to this phenomenon as a ``somme localement finie''. 
In his definition of strict polynomial map, he does not include such an assumption on
bounded degree, but he circumvents it by immediately restricting attention to homogeneous maps. 

We say that a $\NAlg$-polynomial map 
$\phi$ is of \textbf{degree $n$} if $v_{(m_1,\dots,m_k)}=0$ whenever $m_1+\dots+m_k>n$ (independently of the 
linear form  $\sum u_j\otimes x_j$).

\bex 
Let $U=\gen{u_0,u_1, u_2,\dots}$ be free on an infinite basis. The map 
$$ \phi_A\colon U\otimes A \to U\otimes A, \quad \sum u_k\otimes a_k \mapsto \sum u_k\otimes \binom{a_k}{k} $$
is $\NAlg$-polynomial of infinite degree. 
\eex
%

The final theorem of this paper will establish that a map is numerical of degree $n$ if and only if it is 
$\NAlg$-polynomial of degree $n$.

\section{The Universal Numerical Maps} 

When $M$ is a module (or even abelian group), the free module on $M$, as given by 
$$ 
\B[M] = \Set{\sum a_j [x_j] |  a_j\in \B,\ x_j\in M }, 
$$
carries a multiplication $[x][y]=[x+y]$,  
making $\B[M]$ into a commutative, associative algebra with unity $[0]$. 
We consider the $n$'th deviation 
$$
(x_1,\ldots,x_{n+1})\mapsto [x_1\de\cdots\de x_{n+1}].
$$
of the map
$$
M\to \B[M], \qquad x\mapsto [x].
$$ 

\blem 		\label{L: Dev Formula AA} 
In the algebra $\B[M]$, the following formula holds: 
$$
[x_1\de\cdots\de x_{n+1}]  = \left([x_1]-[0]\right) \cdots \left([x_{n+1}]-[0]\right). 
$$
\elem

\bpr 
$$
[x_1\de\cdots\de x_{n+1}] = \sum_{I\subseteq[n+1]} (-1)^{n+1-\abs{I}} \left[\sum_{i\in I} x_i\right] = \left([x_1]-[0]\right) \cdots \left([x_{n+1}]-[0]\right). \qedhere
$$ 
\epr

This leads to a filtration of $\B[M]$, given by the decreasing sequence of ideals
\begin{align}		\label{E: Polynomial modules}
I_n &=  \Ideal{ [x_1\de\cdots\de x_{n+1}] |  x_i\in M } \nonumber \\
& \quad + \Ideal{ [rx] - \sum_{k=0}^n \binom{r}{k}\left[\De_k x\right] | 
 r\in \B, \ x\in M } , \quad n\geq -1. 
\end{align}
Following Passi \cite{PassiI}, we call the quotient algebra 
$\B[M]/I_n$ a \textbf{numerical module}, and denote it by $P_n(M)$. Passi essentially 
considers the case when polynomial 
relations have been factored out (the first summand of \eref{E: Polynomial modules} only, meaning that $I_n$ 
equals the $(n+1)$'st power of the augmentation ideal), and investigates the 
structure of those $P_n(M)$ in \cite{PassiI} and \cite{PassiFunctors}.

\bth 		\label{Th: Universal} 
The map 
$$ 
\delta_n\colon M\to P_n(M), \qquad x\mapsto [x], 
$$
is the universal numerical map of degree $n$, in that every numerical map 
$\phi\colon M\to N$ of degree $n$ has a unique factorisation $\phi=\hat{\phi}\delta_n$ through it: 
$$ 
\xymatrix{
M \ar[r]^-{\delta_n} \ar[dr]_\phi & P_n(M) \dotar[d]^{\hat{\phi}} \\
& N } $$ 
\eth

\bpr 
Given a map $\phi\colon M\to N$, extend it linearly to an homomorphism $\phi\colon \B[M]\to N$.  
The theorem amounts to the trivial observation that $\phi$ is numerical of degree $n$ if and only if it 
annihilates $I_n$. 
\epr

The numerical modules arising from a free module $M$ are given by the next theorem. 

\bth 
In the polynomial algebra $\B[t_1,\ldots,t_k]$, let $J_n$ be the ideal generated by monomials of degree 
greater than $n$. Denote by $(e_i)_{i=1}^k$ 
the canonical basis of $\B^k$. Then 
$$ 
\psi\colon \B[t_1,\ldots,t_k]/J_n \to P_n(\B^k) ,  \qquad
t_1^{m_1}\cdots t_k^{m_k} \mapsto \left[ \De_{m_1} e_1 \de \cdots \de \De_{m_k} e_k \right], 
$$
is an isomorphism of algebras. In particular, $P_n(\B^k)$ is a free module. 
\eth

\bpr 
The map 
$$
\B[t_1,\ldots,t_k] \to P_n(\B^k),   \qquad 
t_1^{m_1}\cdots t_k^{m_k} \mapsto \left[ \De_{m_1} e_1 \de \cdots \de \De_{m_k} e_k \right] = \prod_i ([e_i]-[0])^{m_i}
$$
is clearly an homomorphism of algebras, and, since it annihilates $J_n$, 
it will factor via $\B[t_1,\ldots,t_k]/J_n$. This establishes the existence of $\psi$.

We now define the inverse of $\psi$. 
Each polynomial $f(t)\in\B[t_1,\ldots,t_k]/J_n$ having no constant term is nilpotent, and so the powers 
$$
(1+f)^r = \sum_{m=0}^\infty \binom{r}{m} f^m = \sum_{m=0}^n \binom{r}{m} f^m
$$ 
are defined for any $r\in \B$. Accordingly, 
we may define
\begin{gather*} 
\chi\colon \B[\B^k] \to \B[t_1,\ldots,t_k]/J_n \\
[a_1e_1 + \cdots + a_ke_k]\mapsto (1+t_1)^{a_1}\cdots (1+t_k)^{a_k} + J_n.
\end{gather*}
We write this more succinctly as 
$$ 
[a] \mapsto (1+t)^a + J_n. 
$$
The map $\chi$ is linear by definition, and also multiplicative, since 
$$ 
\chi([a][b]) = \chi([a+b]) = (1+t)^{a+b}  = (1+t)^a(1+t)^b  = \chi([a])\chi([b]).
$$
It maps $I_n$ to $0$, because, when $x_1,\ldots, x_{n+1}\in \B^k$,
\begin{multline*}
\chi([ x_1\de\cdots\de x_{n+1} ]) 
= \sum_{J\subseteq [n+1]} (-1)^{n+1-\abs{J}} \chi\left( \left[ \sum_{j\in J} x_j \right] \right) \\
= \sum_{J\subseteq [n+1]} (-1)^{n+1-\abs{J}} (1+t)^{\sum_{j\in J} x_j}  
= \prod_{j=1}^{n+1} \big( (1+t)^{x_j}-1 \big) = 0,
\end{multline*}
and also, for $r\in \B$ and $x\in \B^k$,
\begin{multline*}
\chi\left([rx] - \sum_{m=0}^n \binom{r}{m}\left[ \De_m x \right]\right) 
= \chi\left([rx] -  \sum_{m=0}^n \binom{r}{m} \sum_{j=0}^m (-1)^{m-j} \binom{m}{j} [jx] \right) \\
= (1+t)^{rx} - \sum_{m=0}^n \binom{r}{m} \sum_{j=0}^m (-1)^{m-j} \binom{m}{j} (1+t)^{jx} \\
= (1+t)^{rx} - \sum_{m=0}^n \binom{r}{m} \big((1+t)^x-1\big)^m.
\end{multline*}
Letting $f=(1+t)^x-1$, the expression equals
$$
(1+f)^{r} - \sum_{m=0}^n \binom{r}{m} f^m = \sum_{m=0}^\infty \binom{r}{m} f^m  - \sum_{m=0}^n \binom{r}{m} f^m,  
$$
which is zero modulo $J_n$.
We therefore have an induced map 
$$\chi\colon  P_n(\B^k) \to \B[t_1,\ldots,t_k]/J_n.$$
The inverse relationship of $\psi$ and $\chi$ is easy to verify. 
\epr


\section{Properties of Numerical Maps} 

Let us elaborate somewhat on the behaviour of numerical maps, and investigate 
one or two of their elementary properties. 
To begin with, we note that the binomial co-efficients themselves, considered as maps $\B\to \B$, are numerical. 

\bth The binomial co-efficient $x\mapsto \binom xn$ is numerical of degree $n$. \eth

\bpr $\binom xn$ is given by a polynomial in the enveloping $\Q$-algebra.  \epr

Next, we prove an alternative characterisation of numericality.

\blem For $r$ in a binomial ring and natural numbers $m\leq n$, the following formula holds:
$$ \sum_{k=m}^n (-1)^k\binom{r}{k}\binom{k}{m} = (-1)^{n}\binom{r}{m}\binom{r-m-1}{n-m}.$$ \label{Lm: Binomial formula} \elem

\bpr 
Induction on $n$. 
\epr

\bth		\label{Th: Poly Map Eq 1} 
Let the map $\phi\colon M\to N$ be polynomial of degree $n$. It is numerical (of degree $n$) if and only if 
it satisfies the equation 
$$ 
\phi(rx) = \sum_{m=0}^n (-1)^{n-m}\binom{r}{m}\binom{r-m-1}{n-m} \phi(mx),\quad r\in \B ,\ x\in M.
$$ 
\eth

\bpr This follows from the lemma:
\begin{align*}
\sum_{k=0}^n \binom{r}{k}\Dev{\phi}{_k x} &= \sum_{k=0}^n \binom{r}{k} \sum_{m=0}^k (-1)^{k-m} \binom{k}{m} \phi(mx) \\
&= \sum_{m=0}^n (-1)^{-m}   \left( \sum_{k=m}^n (-1)^k \binom{r}{k} \binom{k}{m} \right) \phi(mx) \\
&= \sum_{m=0}^n (-1)^{n-m} \binom{r}{m}\binom{r-m-1}{n-m}\phi(mx). \qedhere
\end{align*} 
\epr

The next formula will frequently be found useful. It describes the behaviour of the lower-order 
deviations of a numerical map. 

\bth 		\label{Th: Poly Map Eq 2} 
The map $\phi\colon M\to N$ is numerical of degree $n$ if and only if  
$$ 
\phi(a_1x_1\de\cdots\de a_kx_k) = 
\sum_{\substack{ m_1+\cdots+m_k\leq n \\ m_1,\dots,m_k\geq 1}} \binom{a_1}{m_1}\cdots\binom{a_k}{m_k} 
\phi\left(\De_{m_1} x_1 \de\cdots\de \De_{m_k} x_k\right) 
$$ 
for any $a_i\in\B$ and $x_i\in M$.
\eth

\bpr 
If the equation is satisfied, it follows readily, from a consideration of the two 
special cases when $k=1$ and $a_1=\dots=a_{n+1}=1$, respectively, that the conditions 
for a numerical map are fulfilled. 

Conversely, if $\phi$ is of degree $n$, a calculation in the numerical module $P_n(M)$, 
using Lemma \ref{L: Dev Formula AA}, yields
\begin{align*}
[a_1x_1\de\cdots\de a_kx_k]  &= \left([a_1x_1]-[0]\right) \cdots \left([a_kx_k]-[0]\right) \\
&= \sum_{q_1=1}^n \binom{a_1}{q_1}\left[\De_{q_1} x_1\right] \cdots 
\sum_{q_k=1}^n \binom{a_k}{q_k}\left[\De_{q_k} x_k\right] \\
&= \sum_{q_1=1}^n \cdots \sum_{q_k=1}^n \binom{a_1}{q_1}\cdots \binom{a_k}{q_k} 
([x_1]-[0])^{q_1} \cdots ([x_k]-[0])^{q_k} \\
&= \sum_{q_1=1}^n \cdots \sum_{q_k=1}^n \binom{a_1}{q_1}\cdots \binom{a_k}{q_k} 
\left[\De_{q_1} x_1 \de \cdots \de \De_{q_k} x_k  \right].
\end{align*}
Application of $\phi$ will reach the desired conclusion. 
\epr


The next theorem should be compared with Theorem \ref{S: Roby}, 
which gives the corresponding property for strict polynomial maps. 

\bth 		\label{Th: Char of Poly Maps} 
Let $\phi\colon M\to N$ be numerical of degree $n$ and  
let $u_1, \ldots, u_k \in M$. For any natural numbers $m_1,\dots,m_k$ there is 
an element $v_{(m_1,\dots,m_k)}\in N$, non-zero only if $m_1+\cdots+m_k\leq n$, such that  
\beq 		\label{E: Char of Poly Maps}
\phi(a_1u_1+\cdots+ a_ku_k)= \sum \binom{a_1}{m_1}\cdots \binom{a_k}{m_k} v_{(m_1,\dots,m_k)}
\eeq
for any $a_1,\dots,a_k\in \B$. The elements $v_{(m_1,\dots,m_k)}$ are uniquely defined and, explicitly, 
$$
v_{(m_1,\dots,m_k)} = \phi\left( \De_{m_1} u_1 \de \cdots \de \De_{m_k} u_k\right).
$$

Conversely, the existence of elements $v_{(m_1,\dots,m_k)}$ with the above property 
ensures that $\phi$ is numerical of degree $n$.
\eth

\bpr 
Assume $\phi$ is numerical of degree $n$ and let $I\subseteq [k]$. By the preceding theorem, 
$$
\phi\left( \De_{i\in I} a_iu_i\right) = \sum
\binom{a_1}{m_1} \cdots \binom{a_k}{m_k} \phi\left(\De_{m_1} u_1 \de\cdots\de \De_{m_k} u_k\right),
$$ 
where the sum is taken over all indices $m_1,\dots,m_k$ such that $m_1+\cdots+m_k\leq n$ and 
$m_i\geq 1$ if $i\in I$, but $m_i=0$ if not. 
Hence 
\begin{multline*}
\phi(a_1u_1+\dots+a_ku_k) = \sum_{I\subseteq [k]} \phi\left( \De_{i\in I} a_iu_i\right) \\
= \sum_{m_1+\cdots+m_k\leq n} 
\binom{a_1}{m_1} \cdots \binom{a_k}{m_k} \phi\left(\De_{m_1} u_1 \de\cdots\de \De_{m_k} u_k\right),
\end{multline*}
which shows that the formula \eref{E: Char of Poly Maps} holds as stated. 

Assume now that $\phi$ (not \emph{a priori} numerical) is such that elements $v_{(m_1,\dots,m_k)}$ exist, making 
the formula \eref{E: Char of Poly Maps} valid. 
Fix natural numbers $m_1,\dots,m_k$. By \eref{E: Char of Poly Maps} and Lemma 
\ref{L: Dev Lemma}, we may write 
\begin{multline*}
 \sum_{0\leq s_i\leq m_i} \binom{m_1}{s_1}\cdots \binom{m_k}{s_k} v_{(s_1,\dots,s_k)}  
= \phi(m_1u_1+\dots+m_ku_k)  \\
= \sum_{0\leq s_i\leq m_i} 
\binom{m_1}{s_1}\cdots \binom{m_k}{s_k} \phi\left(\De_{s_1}u_1 \de \cdots \de \De_{s_k} u_k\right).
\end{multline*}
The left-most expression contains the monic term $v_{(m_1,\dots,m_k)}$. Since the remaining indices $(s_1,\dots,s_k)$ 
precede $(m_1,\dots,m_k)$ in the lexicographical ordering on $\N^k$, it follows inductively that 
each 
$$
v_{(m_1,\dots,m_k)} = \phi\left(\De_{m_1}u_1 \de \cdots \de \De_{m_k} u_k\right)
$$
is uniquely determined. Finally, given this explicit expression for $v_{(m_1,\dots,m_k)}$, and the assumption 
that $v_{(m_1,\dots,m_k)}=0$ when $m_1+\cdots+m_k>n$,
it follows readily, by a consideration of the two special cases $k=1$ and $a_1=\cdots=a_{n+1}=1$, respectively, 
that $\phi$ must be numerical of degree $n$.
\epr

\section{Numericality versus $\NAlg$-Polynomiality} 

We shall now tie things together in our main theorem, and 
show that numerical maps and $\NAlg$\hyp polynomial maps are one and the same. 

We refer to Theorem 5 of \cite{BR} for the definition of the free binomial 
algebra $\freenum[\B]{t_1,\dots,t_k}$  on the variables $t_1,\dots,t_k$.

\blem 			\label{L: Zero}
Let $N$ be a $\B$-module and consider a binomial polynomial 
$$
p(t_1,\dots,t_k) = \sum_{m_1,\dots,m_k\in\N} v_{(m_1,\dots,m_k)} \otimes \binom{t_1}{m_1}\cdots\binom{t_k}{m_k} 
\in N\otimes \freenum[\B]{t_1,\dots,t_k}. 
$$
If $p$ vanishes on $\B^k$, then all $v_{(m_1,\dots,m_k)}=0$. 
\elem

\bpr
For a single variable $t$, the polynomial has the form 
$p(t) = \sum_{m} v_{m} \otimes \binom{t}{m}$.
By successively considering
\begin{align*}
0 &= p(0) = v_0 \\
0 &= p(1) = v_0 + v_1 \\
0 &= p(2) = v_0 + 2v_1 + v_2, \quad\text{\&c.,}  
\end{align*}
one finds $v_m=0$ for all $m$. For $k\geq 2$, write 
$$
p(t_1,\dots,t_k) = \sum_{m_k} \left( \sum_{m_1,\dots,m_{k-1}} v_{(m_1,\dots,m_k)} \otimes 
\binom{t_1}{m_1}\cdots\binom{t_{k-1}}{m_{k-1}} \right) \binom{t_k}{m_k} 
$$
and use induction.
\epr

\bth 			\label{Th: Num = NAlg}
The map $\phi\colon M\to N$ is numerical of degree $n$ if and only if 
it may be extended to a (unique) $\NAlg$-polynomial map $M\otimes- \to N\otimes-$ of degree $n$.  
\eth

\bpr 
If $\phi_A\colon M\otimes - \to N\otimes -$ is a $\NAlg$-polynomial map of bounded degree $n$, 
it is clear from the Polynomiality Principle that 
$\phi_\B\colon M\to N$ has the property of Theorem \ref{Th: Char of Poly Maps}, and is therefore numerical. 

We shall show that $\phi_A$ is uniquely specified by the numerical map $\phi_\B$.
Let 
$a_1,\dots,a_k\in \B$ and $u_1,\dots,u_k\in M$. Consider the action of 
$$
\phi\colon M\otimes \freenum[\B]{t_1,\dots,t_k}\to N\otimes\freenum[\B]{t_1,\dots,t_k}.
$$
Write 
$$
\phi\left( \sum u_j\otimes t_j \right) = \sum w_{(m_1,\dots,m_k)} \otimes \binom{t_1}{m_1}\cdots\binom{t_k}{m_k}
$$
for some $w_{(m_1,\dots,m_k)}\in N$. By naturality, $\phi$ commutes with the action of the homomorphism 
$\alpha\colon \freenum[\B]{t_1,\dots,t_k}\to \B$, given by $t_j\mapsto a_j$, 
whence 
\begin{multline*}
\phi_\B\left(\sum a_j u_j\right) 
= \phi_\B\alpha\left( \sum u_j\otimes t_j \right) 
= \alpha\phi\left( \sum u_j\otimes t_j \right) \\
= \alpha\left( \sum w_{(m_1,\dots,m_k)} \otimes \binom{t_1}{m_1}\cdots\binom{t_k}{m_k} \right) 
= \sum  \binom{a_1}{m_1}\cdots\binom{a_k}{m_k} w_{(m_1,\dots,m_k)} .
\end{multline*}
By the uniqueness part of Theorem \ref{Th: Char of Poly Maps}, it must be that 
$$
w_{(m_1,\dots,m_k)}=\phi_\B\left(\De_{m_1} u_1 \de\cdots\de \De_{m_k} u_k\right),
$$ so that necessarily
$$ 
\phi\left( \sum u_j\otimes t_j \right) = \sum \phi_\B\left(\De_{m_1} u_1 \de\cdots\de \De_{m_k} u_k\right)
\otimes \binom{t_1}{m_1}\cdots\binom{t_k}{m_k}.
$$ 
Again by naturality, this same equation must hold for all elements $t_j$ in all binomial algebras, thus 
providing the only possible way to extend $\phi$.

There remains to verify that, given a numerical $\phi_\B$, we may indeed extend it into a 
well-defined natural transformation $\phi_A\colon M\otimes A\to N\otimes A$ by
$$
\phi_A(u_1 \otimes a_1 +\cdots+ u_k \otimes a_k) = \sum  
\phi_\B\left( \De_{m_1} u_1 \de \cdots \de \De_{m_k} u_k\right) \otimes \binom{a_1}{m_1}\cdots \binom{a_k}{m_k}
$$ 
for any $a_1,\dots,a_k$ in any binomial algebra $A$. It is clear this map is natural in $A$. 
To see that the equation defines a well-defined map, consider, for example, the two possible 
images of the tensor $u\otimes a + u\otimes b = u\otimes (a+b)$: 
\begin{gather*}
\phi_A(u \otimes a + u \otimes b) = \sum  
\phi_\B\left( \De_{p} u \de \De_{q} u\right) \otimes \binom{a}{p}\binom{b}{q} \\ 
\phi_A(u \otimes (a + b)) = \sum  
\phi_\B\left( \De_{m} u \right) \otimes \binom{a+b}{m} 
\end{gather*}
The right-hand sides are polynomials in $a,b$, and since they take equal values when $a,b\in\B$ 
by Theorem \ref{Th: Char of Poly Maps} (namely $\phi_\B(au+bu)$), 
they are equal as polynomials by  Lemma~\ref{L: Zero}.
\epr

\end{document}